\documentclass{amsart}
\usepackage{amsmath} 
\usepackage{amssymb}

\newcommand{\forc}{\Vdash} 

\newcommand{\dom}{{\rm Dom}}
\newcommand{\add}{{\rm Add}}

\newcommand{\alp}{{\alpha}}

\newcommand{\kk}{\kappa}

\newcommand{\oo}{{\omega}}
\newcommand{\pb}{\setminus}

\title{A space with only Borel subsets}

\author{Saharon Shelah}
\address{Institute of Mathematics\\
 The Hebrew University of Jerusalem\\
 91904 Jerusalem, Israel\\
 and  Department of Mathematics\\
 Rutgers University\\
 New Brunswick, NJ 08854, USA}
\email{shelah@math.huji.ac.il}
\urladdr{http://www.math.rutgers.edu/$\sim$shelah}
\thanks{The research was partially supported by the Israel Science
Foundation, founded by the Israeli Academy of Sciences and
Humanities. Publication 730}   

\begin{document}
\maketitle

Mikl\'os Laczkovich (Budapest) asked if there exists a Haussdorff (or even
normal) space in which every subset is Borel yet it is not meager.  The
motivation of the last condition is that under ${\rm MA}_\kappa$ every
subspace of the reals of cardinality $\kk$ has the property that all subsets
are ${\rm F}_\sigma$ however Martin's axiom also implies that these subsets
are meager.  Here we answer Laczkovich' question.  I thank Peter Komjath --
the existence of this paper owes much to him.

\bigskip
\noindent{\bf Theorem.} 
         {\sl The following are equiconsistent.}
\begin{enumerate} 
\item {\sl There exists a measurable cardinal.}
\item {\sl There is a non-meager $T_1$ space with no isolated points in 
           which every subset is Borel.}
\item {\sl There is a non-meager $T_4$ space with no isolated points in 
           which every subset is the union of an open and a closed set.}
\end{enumerate}        

\bigskip
\noindent{\bf Proof.} 
Assume first that $\kk$ is measurable in the model $V$. Add $\kk$ Cohen
reals, that is, force with the partial ordering $\add(\oo,\kk)$. Our model
will be $V[G]$ where $G\subseteq\add(\oo,\kk)$ is generic. We first observe
that in $V[G]$ there is a $\kk$-complete ideal on $\kk$ such that the
complete Boolean algebra $P(\kk)/I$ is isomorphic to the Boolean algebra of
the complete closure of $\add(\oo,j(\kk))$ where $j:V\to M$ is the
corresponding elementary embedding. Indeed we let $X\in I$ if and only if
$1\forc \kk\notin j(\tau)$ for some $\tau$ satisfying $X=\tau^G$, that is,
$\tau$ is a name for $X\subseteq\kk$. Moreover, the mapping $X\mapsto
[\![\kk\in j(\tau)]\!]$ is an isomorphism between $P(\kk)/I$
and the regular Boolean algebra of $\add(\oo,j(\kk) \pb \kk)$ (where $\tau$
is a name for $X$). Notice that $|j(\kk)|=2^\kk$.

We observe that this Boolean algebra has the following properties. There
are $2^\kk$ subsets $\{A_\alp:\alp<2^\kk\}$ which are independent mod $I$,
that is, if $s$ is a function from a finite subset of $\kk$ into $\{0,1\}$
then the intersection
\[
 B_s\stackrel{\rm def}{=}\bigcap_{\alp\in\dom(s)} A^{s(\alp)}_\alp 
\]
is not in $I$ (here $A^1=A$ and $A^0=\kk \pb A$). Moreover, if
$A\subseteq\kk$ then there are countably many pairwise contradictory
functions $s_0,s_1,\dots$ as above, such that  
\[ 
 A/I = B_{s_0}/I \lor B_{s_1}/I \lor \dots, 
\] 
that is, $A$ can be written as $B_{s_0}\cup B_{s_1}\cup\dots$
add-and-take-away a set in $I$.
 
By cardinality assumptions we can assume that for every pair $(X,Y)$ of
disjoint members of $I$ there is some $\alp<2^\kk$ with $X\subseteq A_\alp$,
$Y\subseteq \kk \pb A_\alp$.

We define a topology on $\kk$ by declaring the system 
\[ 
 \{ A_\alp \pb Z, A^1 \pb Z:\alp<2^\kk, Z\in I\} 
\] 
a subbasis, or, what is the same, the collection of all sets of the form
$B_s\pb Z$ (where $Z\in I$) a basis.

We prove the following statements on the space.

\bigskip
\noindent{\bf Claim.} 
  {\sl The space has the following properties.}
\begin{enumerate} 
\item {\sl Every set of the form $B_s$ is clopen, every set in $I$ is
closed.}
\item {\sl Every meager set is in $I$.} 
\item {\sl Every set is the union of an open and a closed set. }
\item {\sl The closure of $B_s \pb Z$ is $B_s$.} 
\item {\sl The space is $T_4$.}
\end{enumerate} 

\bigskip
\noindent{\bf Proof.} 
1. Straightforward. 

2. Every set not in $I$ contains a subset of the form $B_s \pb Z$
(by one of the properties of the Boolean algebra mentioned above), which
is open, so every nowhere dense, therefore every meager set is in $I$. 

3. If $A\subseteq\kk$ then $A/I$ can be written as 
$A/I = B_{s_0}/I \lor B_{s_1}/I \lor \dots$ and then clearly 
\[
 A=\Bigl(\bigl(B_{s_0} \pb Z_0 \bigr)\cup 
\bigl(B_{s_1} \pb Z_1 \bigr)\cup \dots\Bigr) \cup Z
\] 
for some sets $Z_0,Z_1,\dots,Z$ in $I$. 
But this is a decomposition into the union of an open and a closed set. 

4. Clear. 

5. Assume we are given the disjoint closed sets $F$ and $F'$. 
They can be written as 
\[
 F = \bigl(B_{s_0} \pb Z_0 \bigr)\cup \bigl(B_{s_1} \pb Z_1 \bigr)\cup
\dots \cup Z
\]
and 
\[
 F' = \bigl(B_{s'_0} \pb Z'_0 \bigr)\cup \bigl(B_{s'_1} \pb Z'_1
\bigr)\cup 
 \dots \cup Z'.
\]
As $F$ and $F'$ are closed, using 4., we can assume that 
\[ 
 Z_0=Z_1=\cdots=Z'_0=Z'_1=\cdots=\emptyset.
\] 
Set $G=B_{s_0}\cup B_{s_1}\cup\dots$, $G'=B_{s'_0}\cup B_{s'_1}\cup\dots$,
then $F=G\cup Z$, $F'=G'\cup Z'$ and these four sets are pairwise disjoint.
It suffices to separate each of the pairs $(G,G')$, $(G,Z')$, $(G',Z)$, and
$(Z,Z')$.  There is no problem with the first case, as $G$, $G'$ are open.
For the last case we use our assumption that some $A_\alp$ separates $Z$ and
$Z'$. For the second, we can assume that $G$ is non empty hence $B_{s_0}$ is
well defined and disjoint to $Z'$, now choose $\alp<\kk$ such that $Z'$ is a
subset of $A_\alp$, and so $G$, $A_\alp \pb B_{s_0}$ is a pair of disjoint
open sets as required. Lastly the third case is similar to the second.

\medskip
We have proved $(1)\longrightarrow(3)$, and $(3)\longrightarrow(2)$ is
trivial; lastly for $(2)\longrightarrow(1)$ assume that $(X,{\mathcal T})$
is a non-meager $T_1$ space with no isolated points in which every subset is
Borel. Let $\{G_\alp:\alp<\tau\}$ be a maximal system of disjoint,
nonempty, meager open sets. Such a system exists by Zorn's lemma. Set
$Y=\bigcup\{G_\alp:\alp<\tau\}$. Clearly, $Y$ is meager. As the boundary of
the open $Y$ is nowhere dense, we get that even the closure of $Y$ is
meager. Then the nonempty subspace $Z=X-\overline{Y}$ has the property that
no nonempty open set is meager and every subset is Borel. If $I$ is the
meager ideal on $Z$ then every subset is equal to some open set mod $I$. We
claim that $I$ is precipitous on $Z$ which implies that in some inner model
there is a measurable cardinal (see \cite{jechprikry},
\cite{jechmagidorprikry}). 

For this, assume that ${\mathcal W}^0, {\mathcal W}^1,\dots$ is a refining
sequence of mod $I$ partitions.  That is, every ${\mathcal W}^n$ is a
maximal system of $I$-almost disjoint open sets, and if $A$ is a member of
some ${\mathcal W}^{n+1}$ then there is some member of ${\mathcal W}^n$
which includes $A$ mod $I$.  We try to find a member $A_n\in{\mathcal W}^n$
such that $\bigcap\{A_n:n<\oo\}$ is nonempty.  To this, observe that the
intersection of two members in ${\mathcal W}^n$ is a meagre open set, hence
is the empty set. Therefore, ${\mathcal W}^n$ is actually a decomposition
of $Z \pb Z_n$ into the union of disjoint open sets where $Z_n$ is a meager
set.  Pick an element in $Z \pb \bigcup\{Z_n:n<\oo\}$ then it is in some
member of ${\mathcal W}^n$ for every $n$ and we are done.


\begin{thebibliography}{9}


\bibitem{jechprikry} T.~Jech, K.~Prikry: 
  {\sl Ideals over uncountable sets: Application of almost disjoint 
  functions and generic ultrapowers}, 
  Memoirs of the A.M.S., {\bf 214}, 1979. 

\bibitem{jechmagidorprikry} T.~Jech, M.~Magidor, W.~Mitchell, K.~Prikry: 
  On precipitous ideals, 
  {\sl Journal of Symbolic Logic} {\bf 45}(1980), 1--8. 

   
\end{thebibliography}
\end{document}